\documentclass[a4paper,12pt]{article}
\usepackage{amsfonts}
\usepackage{amsthm}
\usepackage{amsmath}
\usepackage{amsfonts}
\usepackage{latexsym}
\usepackage{amssymb}

\newtheorem{fed}{Definition}[section]
\newtheorem{teo}[fed]{Theorem}
\newtheorem{lem}[fed]{Lemma}
\newtheorem{cor}[fed]{Corollary}
\newtheorem{pro}[fed]{Proposition}
\theoremstyle{definition}
\newtheorem{rem}[fed]{Remark}
\newtheorem{remdef}[fed]{Remark-Definition}

\def\inc{\subseteq}

\def\beeq{\begin{equation}}
\def\endeq{\end{equation}}
\def\N{\mathbb{N}}
\def\R{\mathbb{R}}
\def\C{\mathbb{C}}

\def\noi{\noindent}
\def\QED{\hfill ${\blacksquare}$}
\def\EOE{\hfill ${\blacktriangle}$}

\def\la{\lambda}
\def\sii{ if and only if }

\def\cD{\mathcal{D}}

\def\ele{\mathcal{L}}

\def\ese{\mathcal{S}}
\def\ete{\mathcal{T}}
\def\eme{\mathcal{M}}
\def\ene{\mathcal{N}}

\def\pqcome{$(P_r,P_l)$-complementable }
\def\pqcom{$(P_r,P_l)$-complementable}

\def\stcome{$(\ese,\ete)$-complementable }
\def\stcom{$(\ese,\ete)$-complementable}

\def\stscom{$(\ese,\ete)$-complementable}
\def\stwcome{$(\ese,\ete)$-weakly complementable }
\def\stwcom{$(\ese,\ete)$-weakly complementable}

\DeclareMathOperator{\difp}{\, \div \, }
\DeclareMathOperator{\sump}{\, : \, }

\newcommand{\pint}[1]{\displaystyle \left \langle #1 \right\rangle}

\newcommand{\hil}{\mathcal{H}}
\newcommand{\op}{L(\mathcal{H})}
\newcommand{\opuno}{L(\mathcal{H}_1)}
\newcommand{\opdos}{L(\mathcal{H}_2)}

\newcommand{\posop}{L(\mathcal{H})^+}
\newcommand{\posopuno}{L(\mathcal{H}_1)^+}
\newcommand{\posopdos}{L(\mathcal{H}_2)^+}

\newcommand{\hiluno}{\mathcal{H}_{1}}
\newcommand{\hildos}{\mathcal{H}_{2}}
\newcommand{\opunodos}{L(\mathcal{H}_{1},\mathcal{H}_{2})}

\newcommand{\conv}{\xrightarrow[n\rightarrow\infty]{}}
\newcommand{\convnorm}{\xrightarrow[n\rightarrow\infty]{\|\,\cdot\,\|}}

\newcommand{\short}[3]{#1_{/(#2 \,,\, #3)}}

\newcommand{\angd}[2]{c_0\left[\,#1,\,#2\,\right]}
\newcommand{\angf}[2]{c\left[\,#1,\,#2\,\right]}

\def\rai{^{1/2}}
\def\inv{^{-1}}
\def\ben{\begin{enumerate}}
\def\een{\end{enumerate}}
\def\barr{\begin{array}}
\def\earr{\end{array}}

\def\bdem{\begin{proof}}
\def\edem{\end{proof} }

\def\prop{Proposition \ref}
\def\teor{Theorem \ref}

\def\coro{Corollary \ref}

\def\defi{Definition \ref}

\newcommand{\ecua}[1]{equation (\ref{#1})}
\newcommand{\peso}[1]{ \quad \text{ \rm  #1 } \quad }

\date{}
\begin{document}


\title{{
{Bilateral Shorted Operators and Parallel Sums}}
}
\author
{Jorge Antezana, Gustavo Corach and Demetrio Stojanoff
\thanks{Partially supported by UBACYT I030, ANPCYT PICT 03-09521, PIP 2188/00 
and UNLP 11/X350.}
\footnote{E-mail addresses: antezana@mate.unlp.edu.ar, gcorach@fi.uba.ar
  and demetrio@mate.unlp.edu.ar}}

\maketitle

\vglue.3truecm

\fontsize {10}{8}\selectfont
\centerline{
{\bf Jorge Antezana and Demetrio Stojanoff }}
 
 \medskip 
\centerline{
Depto. de Matem\'atica, FCE-UNLP,  La Plata, Argentina
and IAM-CONICET.
}

\bigskip

\centerline{{\bf Gustavo Corach }}

\medskip
\centerline{
Depto. de Matem\'atica, FI-UBA  and IAM-CONICET,} \centerline{
Saavedra 15,  Piso 3  (1083), 
Buenos Aires, Argentina.}

\fontsize {12}{14}\selectfont

\vglue.5truecm

\begin{abstract}
In this paper we study shorted operators relative to two 
different  subspaces, for bounded operators on  
infinite dimensional Hilbert spaces. We define two notions of 
``complementability" in the sense of Ando for operators, and study the properties 
of the shorted operators when they can be defined. We use
these facts in order to define and study the notions of parallel sum 
and substraction, in this Hilbertian context. 
\end{abstract}



\section{Introduction}

This paper is devoted to generalize two operations, coming from electrical network theory: parallel sum of matrices and shorting of matrices. In \cite{andersonduf}, W. N. Anderson Jr. and R. J. Duffin defined , for positive (semidefinite) matrices $A$ and $B$ the parallel sum $A\sump B= A(A+B)^\dagger B$. The motivation for studying this operation, and its name, come from the following fact: if two resistive n-port networks, with impedance matrices $A$ and $B$, are connected in parallel, then $A\sump B$ is the impedance matrix of the parallel connection. It should be mentioned that the impedance matrix of a resistive n-port network is a positive (semidefinite) $n\times n$ matrix. On the other side, in \cite{anderson}  Anderson defined, for a positive $n\times n$ matrix $A$ and a subspace $\ese$ of $\C^n$, the shorted matrix of $A$ by $\ese$. 
Just to give an idea about $A_{/\ese}$, suppose that $A$ has the block form $\begin{pmatrix}A_{11}& A_{12}\\ A_{21} & A_{22}\end{pmatrix}$ where $A_{11}$ is a $k\times k$ block and $A_{22}$ is an $(n-k)\times(n-k)$ block. If $\ese$ is the subspace spanned by the first $k$ canonical vectors, then
\[
A_{/\ese}=\begin{pmatrix}
A_{11}-A_{12} A_{22}^{\dagger} A_{21} & 0\\
0 & 0
\end{pmatrix}
\]
where ${\dagger}$ denotes the Moore-Penrose inverse. (Some authors define $A_{/\ese}$ as a linear transformation $\ese\to\ese$ avoiding the zeroes above). The name shorted comes from the fact that it gives the joint impedance of a resistive n-port, some of whose parts have been short circuited. Here $A$ is the impedance matrix of the original network and $A_{/\ese}$ is the impedance matrix of the network after the short circuits. 
Both operations have been studied in Hilbert spaces context (see the historical notes below). 

One of the goals of this paper is to extend the shorting operation to bounded linear operator between two different Hilbert spaces, given a closed subspace on each one. The solution we get, which we call the bilateral shorted operator, comes from a notion of weak complementability, which is a refinement of a finite dimensional notion due to T. Ando \cite{Ando} and generalized by D. Carlson and E. V. Haynworth \cite{[CaHayn]}. The bilateral shorted operator has been studied in finite dimensions by S. K. Mitra and M. L. Puri \cite{[MP]} (see also the papers by H. Goller \cite{Goller} and  Mitra and Prasad \cite{[Mitra1]}, who refined some results of \cite{[MP]}). However, their methods strongly depend on the existence of generalized inverses, so they can not be used for operators with non closed range.

The second goal is to extend parallel summability for two bounded linear operators between different Hilbert spaces. It should be mentioned that C. R. Rao and S. K. Mitra \cite{[MR]},  and Mitra and K. M. Prasad \cite{[Mitra1]} have studied this extension in finite dimensional spaces. Again, generalized inverses are the main tool they use. In order to avoid generalized inverses, we frequently use what we call hereafter Douglas theorem, an extremely useful result due to R. G. Douglas \cite{douglas}, which we describe after fixing some notations.  

In these notes, $\hiluno$ and $\hildos$ denote Hilbert spaces,
$\opunodos$ is the space of all bounded linear operators between
$\hiluno$ and $\hildos$, we write $L(\hil_i)=L(\hil_i,\hil_i)$ and $\posopuno$ 
(resp. $\posopdos$) the cone of all positive operators on $\hiluno$ (resp. $\hildos$.). 
Recall that $C\in\op$ is called positive if $\pint{Cx,\ x}\geq 0$ for every  $ x\in\hil$. 
For every $C\in\opunodos$ its range is denoted by $R(C)$, its
nullspace by $N(C)$. Given two selfadjoint operators $A,B \in \op$, 
 $A\leq B$ means that $B-A\in\posop$ (this is called the usual or L\"owner order).
A projection is an idempotent (bounded linear) operator. Given 
a closed subspace $\ese\subseteq\hil_1$, by $P_\ese \in \opuno$ is denoted
the  orthogonal projection onto $\ese$. 
%
Douglas theorem states that given $A\in \opunodos$ 
and $B \in L( \hil_3 ,\hildos )$, the following conditions are
equivalent:
$$
\mbox{1. $R(B) \subseteq R(A)$,}\ \ \mbox{2. $\exists\ \la\geq 0:$  $BB^* \leq \lambda \ AA^*$}\ \ \mbox{and 3. $\exists\ D \in L( \hil_3 ,\hiluno  ):$  $B = AD$.}
$$
With the additional condition $ R(D) \subseteq N(A)^\perp$, $D$ is unique and 
it is called the \textbf{reduced solution} of the equation 
$AX=B$; it holds that 
$\|D\|^2 = \inf \big\{\lambda \in \mathbb {R}\ : \  BB^* \le \lambda \ AA^* \big\} $ 
and $N(D) = N(B)$. 
 
\smallskip
\noi We shall use the fact that each pair of
closed subspaces $\ese\subseteq\hil_1$ and $\ete\subseteq\hil_2$
induces a representation of elements of $\opunodos$ by $2\times 2$
block matrices. In this sense,  we identify each $A\in\opunodos$ with a
$2\times 2$ matrix, let us say
\begin{equation}\label{matrix}
A =   \begin{array}{c}
  \ete \\ 
  \ete^\bot
\end{array} 
\begin{pmatrix}
  A_{11} & A_{12} \\
  A_{21} & A_{22}
\end{pmatrix}
\begin{array}{c}
  \ese \\
  \ese^\bot
  \end{array},
\end{equation}
where $A_{11}=\left.P_\ete A\right|_\ese \in L(\ese, \ete)$, 
$A_{12}=\left.P_\ete A\right|_{\ese ^\bot} $, 
$A_{21}=\left.P_{\ete^\bot} A\right|_\ese\ $ and
$\ A_{22}=\left.P_{\ete^\bot} A\right|_{\ese^\bot}$.

\bigskip

\noi {\bf Historical survey:} 
In 1947, M.G. Krein \cite{[K]} proved the existence of a maximum (with respect 
to the usual order) of the set 
$
\eme(A,\ese)=\{C\in\posop: C\leq A\ , \ R(C)\subseteq \ese\}.
$ 
Krein used this extremal operator in his theory of extension of symmetric operators. See the paper by Yu. L. Smul'jian \cite{Smul} for more results in similar directions.
Many years later, W. N. Anderson Jr. \cite{anderson} rediscovered, for finite dimensional spaces, the existence of the maximum which will be denoted $A_{/\ese}$ and called the shorted operator of $A$ by $\ese$. 
Some time before, W. N. Anderson and R. J. Duffin \cite{andersonduf} had developed the binary matrix operation called parallel sum: if $A, B\in L(\mathbb{C}^n)^+$ the parallel sum $A\sump B$ is defined by the formula
\[
A\sump B=A(A+B)^\dagger B.
\]
P. Fillmore and J. P. Williams \cite{fill} defined the parallel sum of positive (bounded linear) operators on a Hilbert space $\hil$ and extended many of Anderson-Duffin's results. It should be mentioned that their definition of parallel sum is based on certain Douglas reduced  solutions.
Anderson and G. E. Trapp \cite{andtrapp} defined $A_{/\ese}$ for a positive operator $A$ on $\hil$ and a closed subspace $\ese$ of $\hil$, and proved that $A_{/\ese}$ can be defined by means of parallel sums, and conversely: if $P$ is the orthogonal projection onto $\ese$, then $A\sump nP$ converges to $A_{/\ese}$ in the operator uniform norm; and for $A,B\in\posop$, $A\sump B$ can be defined as the shorted operator of $\begin{pmatrix}A&A\\ A&A+B\end{pmatrix}\in L(\hil\oplus \hil)^+$ by the subspace $\hil\oplus \{0\}$. This is the approach we shall use here.
The shorting of an operator is one of the manifestations of the Schur complement: if $M$ is a square matrix with block form 
$$
M=\begin{pmatrix} A&B\\ C&D\end{pmatrix},
$$
where $A$ and $D$ are also square blocks and $D$ is invertible, the classical Schur complement of $D$ in $M$ is $A-BD^{-1}C$ (see \cite{[Ca]}, \cite{Co} and \cite{[omelet]} for many results, applications and generalizations of this notion). T. Ando \cite{Ando} proposed a generalization of Schur complements which is closer to the idea of the shorted operators. If $A$ is a $n\times n$  complex matrix and $\ese$ is a subspace of $\C^n$, $A$ is called \textbf{$\ese$-complementable} if there are matrices $M_r$ and $M_l$ such that $PM_r=M_r$, $M_lP=M_l$, $PAM_r=PA$ and $M_lAP=AP$. (Here $P$ is the orthogonal projection onto $\ese$). It holds  $AM_r=M_l A M_r= M_l A$ and $AM_r$ does not depends on the particular choice of $M_r$ and $M_l$; Ando calls $A_{\ese}=AM_r$ the Schur compression and $A_{/\ese}=A-A_{\ese}=A-AM_r$ the Schur complement of $A$ with respect to $\ese$. He observes that, if $A$ is a positive $n\times n$ matrix and $\ese$ is the subspace generated by $n-k$ last canonical vectors, then $A_{/\ese}$ has the block form
$$
\begin{pmatrix} A-BD^\dagger C &0\\ 0&0\end{pmatrix},
$$
and therefore, his definition extends the classical Schur complement. D Carlson and E. V. Haynworth \cite{[CaHayn]} observe that a similar construction could be done starting with $A\in C^{n\times m}$ and subspaces $\ese\in\C^n$ and $\ete\in\C^m$. They defined and studied the notion of operators which are complementable with respect to a pair $(\ese,\ete)$.

As Anderson and Duffin remarked in \cite{andersonduf}, the impedance matrix is positive only for resistive networks. In order to study networks with reactive elements, parallel summation and shorting must be extended to not necessarily positive matrices and operators. C. R. Rao and S. K. Mitra \cite{[MR]} defined and studied parallel sums of $m\times n$ matrices and Mitra \cite{[MP]} used their results to define a sort of bilateral shorted operator by two subspaces, one in $\C^n$ and the other in $\C^m$. A common feature in both extension is the use of generalized inverses. It should be mentioned that these constructions can be applied to linear regression problems as in \cite{[MP]},  \cite{[MPut]}, \cite[Appendix]{[Mitra1]}.

\noi We summarize the contents of this paper. In section 3 we study notion of {\it complementability\rm} in  infinite dimensional Hilbert spaces and we define the concept of {\it weakly complementability \rm } (see Definition \ref{definicion de la debil}). 
We also prove in this section the basic  properties of (weakly or not) complementable triples and  we show some criteria for each kind of complementability. In section 4, under some compatibility conditions between the operator $A$ and the subspaces  $\ese$ and $\ete$, we define a bilateral shorted operator $\short{A}{\ese}{\ete}\in\opunodos$, and we study the usual properties of a shorting  operation.
As Mitra \cite{[Mitraminus]} proved for finite dimensional spaces, we show that  
$\short{A}{\ese}{\ete}$ is the maximum of a certain set for a situable order 
(the so called minus order) in $\opunodos$. 
The rest of the paper is devoted the notions of parallel addition and substraction
of operators and their relationship with the shorted operator. The parallel addition is defined by means of the following device, due to Anderson and Trapp: given  $A,B\in\opunodos$, we say that 
$A$ and $B$ are \textbf{weakly parallel summable (resp. parallel summable)}  if 
the triple 
$\begin{pmatrix}   A & A \\   A & A+B \end{pmatrix} \in L\big(\hiluno\oplus \hiluno , 
\hildos \oplus \hildos \big)$, $\hiluno\oplus \{0\}$, $\hildos\oplus \{0\}$  is weakly complementable (resp. complementable). 
In this case we define the \textbf{parallel sum} of $A$ and $B$, denoted by 
$A\sump B\in\opunodos$, as follows: 
\[
\begin{pmatrix}
  A\sump B& 0 \\
  0 & 0
\end{pmatrix}=\left. \begin{pmatrix}
  A & A \\
  A & A+B
\end{pmatrix}\right/ \barr {r}  \\ {(\hiluno\oplus \{0\},\hildos\oplus \{0\})} \earr .
\]
We study the properties of this operator. Again, under the hypothesis of 
summability, all properties of the finite dimensional case are recovered in our context. 
In section 5 we define the notion of parallel substraction, we give some conditions 
which assures its existence and prove some of its properties. 
In section 6, we extend to the bilateral case some well known formulae for the shorted operator in terms of parallel sums and substractions showing that, as for positive operators, parallel and shorting operations can be defined one in terms of the other.  

\section{Preliminaries}

We need the following two definitions of angles between subspaces in a Hilbert space; they are due, respectively, to Friedrichs and Dixmier (see \cite{[Di]} and \cite {[Fr]}, and the excellent survey by Deutsch \cite {[De]}).

\begin{fed}\rm
Given two closed subspaces $\eme$ and $\ene$, the {\it Friedrichs angle}
between $\eme$ and $\ene$ is the angle in $[0,\pi/2]$ whose cosine
is defined by 
\[
\angf{\eme}{\ene}=\sup\Big\{\,|\pint{x , \, y}|:\; x \in
\eme\ominus (\eme\cap \ene), \; y \in \ene\ominus (\eme\cap
\ene)\;\mbox{and}\;\|x\|=\|y\|=1 \Big\}.
\]
The \textit{Dixmier angle} between $\eme$ and $\ene$ is the angle in $[0,\pi/2]$ whose
cosine is defined by
\[
\angd{\eme}{\ene}=\sup\Big\{\,|\pint{x , \, y}|:
\; x\in \eme, \;y\in \ene\;\mbox{and}\;\|x\|=\|y\|=1 \Big\}.
\]
\end{fed}

\noi
The next  proposition collects the results on
angles which are relevant to our work.

\begin{pro}\label{propiedades elementales de los angulos}
\begin{enumerate}
	\item [\rm 1. ] 
	Let $\eme$ and $\ene$ be to closed subspaces of $\hil$. Then
\ben
\item [\rm a. ] 
$\angf{\eme}{\ene}= \angf{\eme^\bot}{\ene^\bot}$
\item [\rm b. ] 
$\angf{\eme}{\ene}<1$ \sii $\eme+\ene$ is closed.
\item [\rm c. ]  
$\hil=\eme^\bot+\ene^\bot$ \sii \ $\angd{\eme}{\ene}<1$.
\een
\item [\rm 2. ] 
(Bouldin \cite{[Bo]}) Given $B\in\opunodos $ and $A\in L(\hildos , \hil_3 )$  with closed range, then
$R(AB)$ is closed \sii $\angf{R(B)}{N(A)}<1 $.
\end{enumerate}
\end{pro}

\section{Complementable operators}

In this section we study complementable operators. We recall different
characterizations of this notion, their extensions to infinite dimensional  Hilbert spaces, and the relationships among them. 
The next definition, due to Carlson and Haynsworth \cite{[CaHayn]}, is an extension of Ando's generalized Schur complement \cite{Ando}.

\begin{fed}\label{definicion de complementable}\rm
Given two projections $P_r\in\opuno$ and $P_l\in\opdos$,
an operator $A\in\opunodos$ is called \pqcome if there exist operators
$M_r\in\opuno$ and $M_l\in\opdos$ such that
\begin{enumerate}
  \item $(I-P_r)M_r=M_r\,$,  \quad $(I-P_l)AM_r=(I-P_l)A$,
  \item $(I-P_l)M_l=M_l$ \quad and \quad  $M_lA(I-P_r)=A(I-P_r)$.
\end{enumerate}
\end{fed}

\medskip

\noi  We shall prove later that this notion only depends on the 
images of $P_r$ and $P_l$.
ike in the finite dimensional case, we have the following alternative 
characterization of complementability. We use freely matrix decompositions 
like \eqref{matrix}. 

\begin{pro}\label{las tres equivalencias}
Let $P_r\in\opuno$ and $P_l\in\opdos$ be two projections whose 
ranges are $\ese$ and $\ete$ respectively. 
Given $A\in\opunodos$, the following statements are equivalent:
\begin{enumerate}
    \item [\rm 1.] $A$ is \pqcom.
    \item [\rm 2.] $R(A_{21})\subseteq R(A_{22})$ and $R(A_{12}^*)\subseteq R(A_{22}^*)$.
  \item [\rm 3.] There exist two projections $\widehat{P}\in \opuno$ and 
  $\widehat{Q}\in \opdos$  such that:
  \begin{align}\label{ayb}
  R(\widehat{P}^{\,*})=\ese&& R(\widehat{Q})=\ete&& R(A\widehat{P})\subseteq \ete&&
  \mbox{and} && R((\widehat{Q}A)^*)\subseteq \ese.
  \end{align}
\end{enumerate}
\end{pro}
\proof 
$1\Rightarrow 2$: By definition \ref{definicion de 
complementable} it holds
$
M_r=
\begin{array}{c}
  \ete \\ 
  \ete^\bot
\end{array} 
\begin{pmatrix}
  \  0 & 0 \ \\
  \ C & D \ 
\end{pmatrix} \begin{array}{c}
  \ese \\
  \ese^\bot
  \end{array},
$ 
and $A_{21}=A_{22}C$. 
Hence $R(A_{21})\subseteq R(A_{22})\,$. Similar arguments show that 
$R(A_{12}^*)\subseteq R(A_{22}^*)$.

\noi 
$2\Rightarrow 3$: Let $E$ and $F$ be the reduced solutions of $A_{21}=A_{22}X$ 
and $A_{12}^*=A_{22}^*X$, respectively. Note that $E \in L(\ese , \ese ^\perp )$ 
and  $F \in L(\ete^\perp  , \ete )$. If
\begin{align*}
\widehat{P}=\begin{pmatrix}
   I & 0 \\
  -E & 0
\end{pmatrix}\begin{array}{c}
  \ese \\
  \ese^\bot
  \end{array} \in \opuno 
&& \mbox{and} && \widehat{Q}=\begin{pmatrix}
   I & -F \\
   0 & 0
\end{pmatrix}\begin{array}{c}
  \ete \\
  \ete^\bot
  \end{array} \in \opdos \ ,
\end{align*}
easy computations show	 that these projections satisfy \ecua{ayb}.

\noi $3\Rightarrow 1$: Define $M_r=I-\widehat{P}$ and
$M_l=I-\widehat{Q}^*$.
Then $R(M_r)=\ese^\bot$ and $R(M_l)=\ete^\bot$, so conditions 1. and
3. of Definition \ref{definicion de complementable} are satisfied. 
On the other hand 
\begin{align*}
(I-Q)AM_r=(I-Q)A(I-\widehat{P})=(I-Q)A-(I-Q)A\widehat{P}=(I-Q)A, \peso{and}
\end{align*}
\begin{align*}
M_lA(I-P)=(I-\widehat{Q}^*)A(I-P)=A(I-P)-\widehat{Q}^*A(I-P)=A(I-P).
\end{align*}
This shows that conditions 2. and 4. of Definition
\ref{definicion de complementable} also hold. \QED

\medskip

\noi The next characterization has been considered in \cite{CMS2}
for self-adjoint operators in a Hilbert space. We prove an
extension to our general setting.
\begin{pro}\label{caracterizacion de la fuerte con subespacios}
Let $P_r\in\opuno$ and $P_l\in\opdos$ be two projections with
ranges $\ese$ and $\ete$, respectively, and let
$A\in\opunodos$. Then the following statements are equivalent:
\begin{enumerate}
  \item [\rm 1.] $A$ is \pqcom.
  \item [\rm 2.] $\hiluno=\ese^\bot+A^{-1}(\ete)$ and
  $\hildos=\ete^\bot+A^{*-1}(\ese)$.
  \item [\rm 3.] $\angd{\ese}{\overline{A^*(\ete^\bot)}}<1$ and  $\angd{\ete}{\overline{A(\ese^\bot)}}<1$.
\end{enumerate}
\end{pro}
\proof 
$1\Longleftrightarrow 2$: Suppose that $A$ is \pqcom. By Proposition
\ref{las tres equivalencias}, there exists a projection $P$ such
that $R(P^*)=\ese$ and $R(AP)\subseteq \ete$. Then,
$N(P)=\ese^\bot$ and $R(P)\subseteq A^{-1}(\ete)$. Hence
$\hiluno=\ese^\bot+A^{-1}(\ete)$.

\noi Conversely, suppose that $\hiluno=\ese^\bot+A^{-1}(\ete)$
and define $\ene=\ese^\bot\cap A^{-1}(\ete)$. Then
$\hiluno=\ese^\bot\oplus ( A^{-1}(\ete)\ominus \ene)$. Let
$\widehat{P}$ be the oblique projection onto $A^{-1}(\ete)\ominus
\ene$ parallel to  $\ese^\bot$. Then, $
R(\widehat{P}^{\,*})=N(\widehat{P})^\bot=R(I-\widehat{P})^\bot=\ese
$, and $R(A\widehat{P})\subseteq \ete$ because
$R(\widehat{P})=A^{-1}(\ete)\ominus \ene$.
Similar arguments show that the existence of a projection $Q$
such that $R(Q)=\ete$ and $R((QA)^*)\subseteq\ese$ is equivalent
to the identity $\hildos=\ete^\bot+A^{*-1}(\ese)$.

\noi 
$2\Longleftrightarrow 3$: It follows from Proposition \ref{propiedades elementales de los angulos} (item 3) and 
the equality $A^*(\ete^\bot)^\bot = A\inv(\ete )$. 
\QED

\begin{remdef}\rm 
\prop{caracterizacion de la fuerte con subespacios}, 
as well as  \prop{las tres equivalencias}, 
shows that the notion of \pqcome operators only depends on $R(P_l)$ and $R(P_r)$. 
Hence, from now on we shall say that an operator $A\in\opunodos$ is \textbf{\stcome} 
instead of \pqcom. \EOE
 \end{remdef}

\noi In finite dimensional spaces, given a fixed subspace $\ese$, every positive 
operator $A$ is $(\ese,\, \ese)$-complementable. Indeed, if
$ 
A=\begin{pmatrix}
  A_{11} & A_{12} \\
  A_{21} & A_{22}
\end{pmatrix}
\begin{array}{c}
  \ese \\
  \ese^\bot
\end{array},
$ 
the inclusion $R(A_{21})\subseteq R(A_{22})$ always holds (see \cite{Smul} for details). However, in infinite dimensional Hilbert spaces, only the inclusion $R(A_{21})\subseteq R(A_{22}^{1/2})$ holds in general. As $R(A_{22})= R(A_{22}^{1/2})$ if and only if $A_{22}$ has closed range (which is the case in finite dimensional spaces), it is not difficult to 
find examples of positive operators which are not $(\ese,\, \ese)$-complementable (e.g., see example 5.5 of \cite{CMS2}). 
For this reason we consider the following weaker notion of complementability:

\begin{fed}\label{definicion de la debil}\rm
Let $\ese\subseteq\hiluno$ and $\ete\subseteq\hildos$  be closed subspaces.
An operator $A\in\opunodos$ is called {\bf \stwcome} if 
$$
R(A_{21})\subseteq R(|A_{22}^*|^{1/2}) \peso{ and } 
R(A_{12}^*)\subseteq R(|A_{22}|^{1/2})\ ,$$
(according to the matrix decomposition of $A$ given in equation \eqref{matrix}). 
\end{fed}

\begin{rem}
Observe that, by Douglas theorem, $R(|A_{22}^*|) = R(A_{22}) $ and 
$R(|A_{22}|) = R(A_{22}^*) $. Therefore this notion is, indeed, weaker 
than the previously defined notion of complementability. However, if $R(A_{22})$ is closed, then  $R(|A_{22}^*|)$ is also closed and $R(A_{22})=R(|A_{22}^*|)=R(|A_{22}^*|^{1/2})$. Thus, both notions of complementability coincide. 
\EOE
\end{rem}

\noi As an easy consequence of Douglas theorem,we get the next alternative 
characterizations of \stwcome operators.

\begin{pro}\label{caracterizaciones de la debil}
Given $A\in\opunodos$, and closed subspaces $\ese\subseteq\hiluno$
,$\ete\subseteq\hildos$, then the following
statements are equivalent:
\ben 
  \item[\rm 1.] $A$ is \stwcom.
  \item[\rm 2.] If $A_{22}=U|A_{22}|$ is the polar decomposition of $A_{22}$, then the equations $A_{21}=|A_{22}^*|^{1/2}U X$ and $A_{12}^*=|A_{22}|^{1/2}Y$ have solutions.
  \item[\rm 3.]$\displaystyle
  \sup_{x\in\ese}\frac{\|A_{21}x\|^2}{\pint{|A_{22}^*|x,\,x}}<\infty\;$ and
           $\displaystyle
  \;\sup_{y\in\ete}\frac{\|A_{12}^*y\|^2}{\pint{|A_{22}|y,\,y}}<\infty$. 
\een 
\end{pro}

\section{Shorted Operators}


Recall that, in the classic case, i.e., if $\hiluno = \hildos = \hil$, 
$\ese = \ete$ and $A \in \posop$, Anderson and Trapp \cite{andtrapp}
proved that $A_{/\ese}=\left(\begin{matrix}
A_{11}-C^*C & 0\\
     0      & 0
\end{matrix}\right)$, where  $C$ is the reduced solution of
$A_{22}^{1/2}X=A_{21}$.
Following this approach, we shall extend the notion of shorted
operators to operators between 
Hilbert spaces $\hiluno$ and $\hildos\,$. 
Throughout this section, $\ese \subseteq \hiluno$ and $\ete\subseteq\hildos$
are closed subspaces and each operator $A\in\opunodos$ is
identified with a $2\times 2$ matrix induced by these subspaces, 
as in \eqref{matrix}.

\medskip
\begin{fed}\rm
Let $A\in\opunodos$ be \stwcom, and let $F$ and $E$ be the reduced solutions of
the equations $A_{21}=|A_{22}^*|^{1/2}U X$ and $A_{12}^*=|A_{22}|^{1/2}X$,
respectively, where $U$ is the partial isometry of the polar decomposition
of $A_{22}$. The \textbf{bilateral shorted operator} 
of $A$ to the subspaces $\ese$ and $\ete$ is
\[
\short{A}{\ese}{\ete}=\begin{pmatrix}
A_{11}-F^*E & 0\\
     0      & 0
\end{pmatrix}.
\]
\end{fed}

\begin{rem}
If $A_{22}$ has closed range, then $\short{A}{\ese}{\ete}=\begin{pmatrix}
  A_{11}-A_{12}A_{22}^{\dagger}A_{21} & 0\\
     0      & 0
\end{pmatrix}$. \EOE
\end{rem}

\noi In the following proposition we collect some basic properties of shorted
operators. The proof is straightforward.

\begin{pro}\label{propiedades directas}
Let $A\in\opunodos$ be \stwcom. Then
\ben
  \item[\rm 1.] for every $\alpha\in\mathbb{C}$, $\alpha A$ is \stwcom, and
  $\short{(\alpha A)}{\ese}{\ete}=\alpha(\short{A}{\ese}{\ete})$.
  \item[\rm 2.] $A^*$ is $(\ete,\ese)$- weakly complementable, and
  $(\short{A}{\ese}{\ete})^*=\short{(A^*)}{\ete}{\ese}$.
  \item[\rm 3.] $\short{A}{\ese}{\ete}$ is \stwcome and
  $\short{(\short{A}{\ese}{\ete})}{\ese}{\ete}=\short{A}{\ese}{\ete}$.
  \item[\rm 4.] if $A=A^*$ and $\ese=\ete$, then
  $\short{A}{\ese}{\ese}$ is self-adjoint.
\een
\end{pro}

\noi The next Proposition is similar to Theorem 1 by Butler and Morley in \cite{butmor}. For the reader's convenience, we include a proof adapted to our setting. First we need a lemma.

\begin{lem}\label{el lema ese}
Let $A,B\in\opunodos$ be such that $R(A^*)\subseteq R(|B|^{1/2})$.
Suppose that there exist a sequence $\{y_n\}$ in $\hiluno$,
$d\in\hildos$ and a positive number $M$ satisfying
$$
A y_n\xrightarrow[n\rightarrow\infty]{}d\ ,\quad B
y_n\xrightarrow[n\rightarrow\infty]{}0\
,\peso{and}\pint{|B|y_n,y_n}\leq M.
$$
Then $d=0$.
\end{lem}
\bdem
Since $\||B|^{1/2}y_n\|^2=\pint{|B|y_n,y_n}\leq M$, we can suppose, 
with no loss of generality, that there exists
$z\in\hil$ such that $|B|^{1/2}y_n\conv z$. As $By_n\conv
0$, it holds $z\in N(|B|^{1/2})$. Let $C$ be the reduced
solution of $A^*=|B|^{1/2}X$. Since $Ay_n\conv d$ we get
$C^*z=d$ and $d=0$ because $N(|B^{1/2}|)\subseteq N(C^*)$.\edem

\begin{pro}\label{teorema uno}
Let $A\in\opunodos$ be \stwcom. Then, given $x\in\ese$ there exist a
sequence $\{y_n\}\subseteq \ese^\bot$ and a positive number
$M$ such that
$$
A\begin{pmatrix}
x\\
y_n
\end{pmatrix}\xrightarrow[n\rightarrow\infty]{}
\short{A}{\ese}{\ete}
\begin{pmatrix}
x\\ 0
\end{pmatrix}
\ ,\peso{and}\pint{|A_{22}|y_n,\;y_n}\leq M \ , \ \ n\in \N .
$$
Conversely, if there exists a sequence $\{z_n\}$ in $\ese^\bot$,
$d\in\ete$, and a positive number $M$ such that
\begin{equation}\label{nuevazo2}
A\begin{pmatrix}
x\\
z_n
\end{pmatrix}\xrightarrow[n\rightarrow\infty]{}
\begin{pmatrix}
d\\
0
\end{pmatrix}
\ ,\peso{and}\pint{|A_{22}|z_n,\;z_n}\leq M,
\end{equation}

\noi then $\begin{pmatrix}
d\\
0
\end{pmatrix}
=\short{A}{\ese}{\ete}\begin{pmatrix}
x\\
0
\end{pmatrix}
$.
\end{pro}
\bdem
Let $E$ and $F$ be the reduced solutions of
$A_{21}=|A_{22}^*|^{1/2}UX$ and $A_{12}^*=|A_{22}|^{1/2}X\ ,
$ respectively. As $
R(E)\subseteq \overline{R(U^*|A_{22}^*|^{1/2})}
=\overline{R(|A_{22}|^{1/2})}$, given $x\in\hiluno$ there is
a sequence $\{y_n\}$ such that $|A_{22}|^{1/2}y_n\conv -Ex$. Then
\begin{eqnarray*}
A_{21}x+A_{22} y_n&=& 
A_{21}x+U\,|A_{22}|^{1/2}|A_{22}|^{1/2}y_n \\ & = &  A_{21}x+|A_{22}^*|^{1/2}U(|A_{22}|^{1/2}y_n)
\xrightarrow[n\rightarrow\infty]{}0,  \peso{and}\\
A_{11}x+A_{12} y_n&=& A_{11}x+F^*|A_{22}|^{1/2}y_n \\& = &  A_{21}x+F^*(|A_{22}|^{1/2}y_n)\xrightarrow[n\rightarrow\infty]{}
\short{A}{\ese}{\ete}(x).
\end{eqnarray*}

\noi Finally, since the sequence $\{|A_{22}|^{1/2}y_n\}$ converges, 
then $\sup_{n\in\mathbb{N}}\pint{|A_{22}|y_n,\,y_n} <\infty$.
Converselly, suppose that there exists another sequence $\{z_n\}$
in $\ese^\bot$ which satisfies \eqref{nuevazo2}. If $w_n=y_n-z_n$,
then $\pint{|A_{22}|w_n,\,w_n}\leq K$. On the other hand,
$A_{11}\, x+A_{12} \, y_n \conv \ d$ and
$A_{11}\, x+A_{12} \, z_n \conv \ \short{A}{\ese}{\ete}
$.
Therefore, $A_{12}\, w_n \conv d-\short{A}{\ese}{\ete}(x)$. 
In a similar way, we obtain that
$A_{22}\, w_n \conv 0$. Therefore,
by Lemma \ref{el lema ese}, we get that
$d=\short{A}{\ese}{\ete}(x)$. \edem

\begin{cor}\label{rango en el caso weak}
Let $A\in\opunodos$ be \stwcom. Then
\begin{eqnarray}
R(A)\cap\ete\subseteq &R(\short{A}{\ese}{\ete})&\subseteq
\overline{R(A)}\cap\ete \label{inclusion uno}\\
R(A^*)\cap\ese\subseteq &R((\short{A}{\ese}{\ete})^*)&\subseteq
\overline{R(A^*)}\cap\ese \label{inclusion dos}
\end{eqnarray}
In particular, $R(\short{A}{\ese}{\ete})=R(A)\cap \ete$ and $R((\short{A}{\ese}{\ete})^*)=R(A^*)\cap\ese$ if $R(A)$ is closed.  
\end{cor}
\bdem
Firstly, we shall prove that $R(\short{A}{\ese}{\ete})\subseteq
\overline{R(A)}\cap\ete$. Clearly, by definition,
$R(\short{A}{\ese}{\ete})$ $\subseteq \ete$. On the other hand, given
$x\in\hiluno$, by Proposition \ref{teorema uno}, there exists a
sequence $\{y_n\}$ in $\ese^\bot$ such that
$
A\begin{pmatrix}
   Px \\
   y_n \
\end{pmatrix}\xrightarrow[n\rightarrow\infty]{}
\short{A}{\ese}{\ete}
\begin{pmatrix}
   x \\
   0 \\
\end{pmatrix}.
$ 
Thus $R(\short{A}{\ese}{\ete})\subseteq \overline{R(A)}$.

\medskip

\noi In order to prove the first inclusion in (\ref{inclusion uno}),
take $x\in R(A)\cap\ete$, and let $z\in\hiluno$ such that $Az=x$.
If $P$ is the orthogonal projection onto $\ese$, then 
$
A
\begin{pmatrix}
   Pz \\
   z-Pz \
   \end{pmatrix}=\begin{pmatrix}
       x \\
       0 \\
\end{pmatrix},
$ 
and, by Proposition \ref{teorema uno}, we get 
$\short{A}{\ese}{\ete}(Pz)=x$. The other inclusions follow in the
same way. \edem

\noi Next, we shall study the shorting operation on \stcome operators. 

\begin{pro}\label{teorema uno para compatibles}
Let $A\in\opunodos$ be \stscom. For every
$x\in\ese$ there exists $y \in\ese^\bot$ such that
$$
A\begin{pmatrix}
  x \\
  y
\end{pmatrix}=\short{A}{\ese}{\ete}\begin{pmatrix}
  x \\
  0
\end{pmatrix}.
$$
Moreover, there exist projections $P\in\opuno$ and $Q\in\opdos$
such that
\begin{equation}\label{QA=AP}
R(P^*)=\ese \ ,  \quad R(Q)=\ete \peso{ and }
QA=AP=\short{A}{\ese}{\ete} \ .
\end{equation}
\end{pro}
\bdem
By Proposition \ref{las tres equivalencias},
there exists a projector $P\in\opuno$ such that $R(P^*)=\ese$ and
$R(AP)\subseteq \ete$. The matrix decomposition of
$P$ with respect to $\ese$ is
$\begin{pmatrix}
  I & 0 \\
  E & 0
\end{pmatrix},
$ where  $I$ is the identity operator of $\ese$ and $E\in L(\ese,\ese^\bot)$.
If $x\in \ese$ and $y=Ex$, then $A\begin{pmatrix}
  x \\
  y
\end{pmatrix}=AP\begin{pmatrix}
  x \\
  0
\end{pmatrix}\in\ete$. If $z_n=y$ for every $n\in\N$, the sequence $\{z_n\}$ satisfies \eqref{nuevazo2}. Hence, by Proposition \ref{teorema uno}, $A\begin{pmatrix}
  x \\
  y
\end{pmatrix}=\short{A}{\ese}{\ete}\begin{pmatrix}
  x \\
  0
\end{pmatrix}$. Therefore $AP=\short{A}{\ese}{\ete}$. In a similar way it can be proved that there exists $Q\in\opdos$
with $R(Q)=\ete$ such that $QA=\short{A}{\ese}{\ete}$ \edem

\begin{rem}\label{teorema dos para compatibles}
Note that we actually prove that if there exists a projection $P$
such that $R(P^*)=\ese$ and $R(AP)\subseteq \ete$, then, by
Proposition \ref{teorema uno}, $AP=\short{A}{\ese}{\ete}$. This
result, for positive operators, appeared in \cite{CMS2},
where the role of $P$ is played by a so-called $A$-selfadjoint
projection, i.e., a projection which is selfadjoint with respect
to the sesquilinear form $\pint{x,\ y}_A=\pint{Ax,\ y}$. The
reader is referred to \cite{[CMS2]},\cite{[HN]} for more information 
about $A$-selfadjoint projections.\EOE
\end{rem}

\begin{cor}\label{Rango del shorted cuando es compatible.}
Let $A\in\opunodos$ be \stscom. Then,
$$
R(\short{A}{\ese}{\ete})  =R(A)\cap\ete \quad \text{and}  \quad
N(\short{A}{\ese}{\ete})=\ese^\bot+N(A) .
$$
\end{cor}
\bdem
By \coro{rango en el caso weak}, it holds
$R(A)\cap\ete \subseteq R(\short{A}{\ese}{\ete})$ and
$$
\ese^\bot+N(A)\subseteq   \big(\ese \cap \overline{R(A^*)}\ \big) ^\perp
\subseteq R(\short{A}{\ese}{\ete}\ ^*)^\perp =   N(\short{A}{\ese}{\ete}) .
$$
On the other hand, by \prop{teorema uno para compatibles}, there
exist two projections $P\in\opuno$ and $Q\in\opdos$ which
satisfy \ecua{QA=AP}. Hence, $R(\short{A}{\ese}{\ete}) = R(AP) \subseteq R(A)$, and
\[
N(\short{A}{\ese}{\ete})=N(AP)=N(P)\oplus \Big(
R(P)\cap N(A) \Big) \subseteq \ese ^\perp+N(A),
\]
because $N(P)=R(P^*)^\bot=\ese^\bot$. 
\edem
\begin{rem}
If $A\in\opunodos$ be \stcom, then, by  \coro{Rango del shorted cuando es compatible.}, 
the subspaces $\ese^\bot+N(A)$, $\ \ese+N(A)^\bot$, $\ \ete^\bot + R(A)^\perp$ and  
$\ \ete + \overline{R(A)}$ must be closed. Moreover, if $R(A)$ is closed then, 
by \prop{propiedades elementales de los angulos}, 
$A( \ese^\perp ) $, $A^*(\ete^\perp ) $, and $R(A_{22})$ are also closed.
Hence, in this case, generalized inverse methods can be used. 
Nevertheless, by using the approach developed in this work,
one can get almost all known properties of the Schur complements
in  finite dimensional spaces, for complementable operators
in general Hilbert spaces, including those operators whose
ranges are not closed.
\EOE
\end{rem}

\bigskip
\noi {\bf The minus partial order.} In \cite{[Mitraminus]}, Mitra proved 
(for matrices in $\C^{m\times n}$) 
that $\short{A}{\ese}{\ete}$ is the unique maximum of the set
\[
\eme^-(A,\ese,\ete)=\Big\{C\in\C^{m\times n}: \ C\leq^-A,\ \ R(C)\subseteq \ete
\peso{and} R(C^*)\subseteq \ese\Big\},
\]
where the partial ordering is the so called minus order:  $C\leq^- A$ if 
$$
R(C)\cap R(A-C)=\{0\}\peso{and} R(C^*)\cap R(A^*-C^*)=\{0\}.
$$ 
A similar result can be obtained in our setting with suitable changes. Firstly, we need to extend the minus order to infinite dimensional Hilbert spaces:

\begin{fed}\label{minus order} \rm
Given $A, B\in\opunodos$, we write $A\leq^- B$ if:
\begin{align*}
	\mbox{\text{(a) }}\ \, \angd{\overline{R(A)}}{\overline{R(B-A)}}<1&& \mbox{and} &&
	\mbox{\text{(b) }}\ \, \angd{\overline{R(A^*)}}{\overline{R(B^*-A^*)}}<1 \ .
\end{align*}
\end{fed}

\begin{rem}
In the finite dimensional case, condition (a) is equivalent to $R(A)\cap R(B-A)=\{0\}$ and condition (b) is equivalent to $R(A^*)\cap R(B^*-A^*)=\{0\}$. So, Definition \ref{minus order} extends   the (finite dimensional) minus order. Also notice that $A\leq^- B$ if and only if $A^*\leq^- B^*$, by the symmetry of conditions (a) and (b). 
\EOE
\end{rem}

\noi The next Proposition provides equivalent conditions to condition (a) in 
\defi{minus order}, which are simpler to handle. 
A similar result for condition (b) can be obtained by taking adjoints.

\begin{pro}\label{equivalencias para el minus}
Given $A,B\in\opunodos$, the following statements are equivalent:
\begin{enumerate}
	\item $\angd{\overline{R(A)}}{\overline{R(B-A)}}<1$.
	\item There exists a projection $Q\in \opdos$ such that $R(Q)=\overline{R(A)}$  and  $A=QB$.
	\item There exists a projection $Q\in \opdos$ such that  $A=QB$.
\end{enumerate}
\end{pro}
\proof $1\Longrightarrow 2:\ $ Let $\ele=\overline{R(A)}\oplus\overline{R(B-A)}$, which is closed by \prop{propiedades elementales de los angulos}. 
	Let $Q \in \opdos$ be the projection with  $R(Q)= \overline{R(A)}$ and $N(Q) 
	= \overline{R(B-A)}\oplus \ele^\perp$. Then, $QB=Q\Big((B-A)+A\Big)=QA=A$.
	
	\noi $2\Longrightarrow 3:\ $ It is apparent.
	
	\noi $3\Longrightarrow 1:\ $ Since $A=QB$ and $B-A=(I-Q)B$, it holds that 
	$R(A)\subseteq R(Q)$ and $R(B-A)\subseteq R(I-Q)=N(Q)$. Hence, $\angd{\overline{R(A)}}{\overline{R(B-A)}}\le \angd{R(Q)}{N(Q)}<1$. \QED

\begin{cor}\label{ring}
Let $A, B \in \opunodos$. 
\ben
\item If $A\leq^- B$, then $R(A) \inc R(B)$ and $R(A^*) \inc R(B^*)$.
\item The relation $\leq^-$ is a partial order (i.e. it is reflexive, 
antisymmetric and transitive). 
\item If $A\leq^-B$ and $B$ is a projection, then $A$ is also a projection.
\een
\end{cor}
\bdem
The first two statements follow easily from Proposition \ref{equivalencias para el minus}. 
If $A\leq^-B$ and $B^2 = B$, by \prop{equivalencias para el minus} 
applied to $A$ and $B$ (resp $A^*$ and $B^*$) 
there exist projections $P$ and $Q$ such that $R(P^*)=\overline{R(A^*)}$, $R(Q)=\overline{R(A)}$ and $A=QB=BP$. Then
$A^2=(QB)(BP)=QBP=A$.
\edem

\begin{teo}\label{mitraing}
Let $A\in \opunodos$ be \stcom, and let
\[
\eme^-(A,\ese,\ete)=\Big\{C\in\opunodos: \ C\leq^-A,\ \ R(C)\subseteq \ete
\peso{and} R(C^*)\subseteq \ese\Big\}.
\]
Then, $\displaystyle \short{A}{\ese}{\ete}=\max_{\leq^-}\ \eme^-(A,\ese,\ete)$.
\end{teo}
\bdem
By Propositions \ref{teorema uno para compatibles} and  \ref{equivalencias para el minus}, 
we know that  $\short{A}{\ese}{\ete}\leq^- A$. On the other hand, 
by Corollary \ref{Rango del shorted cuando es compatible.}, 
$R(\short{A}{\ese}{\ete})\subseteq \ete$ and 
$R((\short{A}{\ese}{\ete})^*)\subseteq \ese$. 
Hence, $\short{A}{\ese}{\ete}\in \eme^-(A,\ese,\ete)$. 
On the other hand, given $C\in \eme^-(A,\ese,\ete)$, there exists a projection 
$E\in\opdos$ such that $C=EA$. Let $P\in\opuno$ be a projection as in 
Proposition \ref{teorema uno para compatibles} such that $R(P^*)=\ese$ 
and $\short{A}{\ese}{\ete}=AP$. The inclusion $R(C^*)\subseteq \ese$ 
implies that $P^*C^*=C^*$. Therefore
\[
 C=CP=EAP=E\short{A}{\ese}{\ete}.
\]
In a similar way, there exists a projection $F$ such that $C^*=F(\short{A}{\ese}{\ete})^*$. So, by Proposition \ref{equivalencias para el minus}, $C\leq^- \short{A}{\ese}{\ete}$.
\edem
\begin{cor}
Let $A \in\op$ be a projection.
If $\ese, \ete \inc \hil$ are closed subspaces such that
$A$ is \stcom , then $N(A) + \ese^\perp$ is closed, 
$$
\hil =  \Big(R(A) \cap \ete\Big) \oplus \Big(N(A) + \ese^\perp \Big) \ ,
$$
and $\short{A}{\ese}{\ete} $ is the  projection given by this decomposition.  
\end{cor}
\bdem 
By \teor {mitraing}, $\short{A}{\ese}{\ete}  \le^- A$. 
Hence it must be a projection by \coro{ring}.
The rest  of the statement follows from \coro{Rango del shorted cuando es compatible.}.
 \edem

\medskip
\noi Next, we shall study the effect of shorting a shorted
operator. The following proposition was proved for selfadjoint
operators by Ando (see \cite{Ando}).

\begin{cor}\label{Teorema sobre el shorted de un shorted}
Let $A\in\opunodos$, and consider closed subspaces $\ese,\widehat\ese$ of $\hil_1$
and $\ete,\widehat\ete$ of $\hil_2$. Then,
it holds
\begin{equation}\label{shorteds iterados}
\short{\left(\short{A}{\ese}{\ete}\right)}{\widehat\ese}{\widehat\ete}=
\short{A}{\ese\cap\widehat\ese}{\ete\cap\widehat\ete}.
\end{equation}
if every operator is complementable with respect to the corresponding pair of subspaces.
\end{cor}
\bdem
Strightforward calculations show that $\eme^-(\short{A}{\ese}{\ete},\widehat\ese,\widehat\ete)=\eme^-(A,\ese\cap\widehat\ese,\ete\cap\widehat\ete)$. Then apply \teor{mitraing}. 
\edem

\begin{rem}
Actually, the last result holds with weaker hypothesis;
in fact, it is only needed that any two of the three shorted
operators exist. 
The reader is
referred to \cite{Ando} for the proofs of these facts. Ando's proof,
valid for a single subspace ($\ese=\ete$), can be easily extended
to our setting. \EOE
\end{rem}

\section{Parallel sum and parallel substraction}
The device of parallel sum of matrices has been
developed by Anderson and Duffin in \cite{andersonduf}. The
extension to general Hilbert spaces is due to
Anderson and Trapp in \cite{andtrapp} (see also \cite{[Mitra1]} and \cite{[Mitra2]}). 
The key idea was to define
parallel sum through shorted operators.
In this section, we shall define parallel sum between operators following the ideas of Anderson and Trapp (see, in particular, \cite{andtrapp} section 4). Even in the scalar case, not every two operators
are summable. So, we need to define the concept of summable operators.

\begin{fed}\rm
Let $A,B\in\opunodos$. We say that $A$ and $B$ are \textbf{weakly parallel
summable} if the next range inclusions hold:
\ben
\item [1. ] $R(A)\subseteq R(|A^*+B^*|^{1/2})$ and $R(B)\subseteq R(|A^*+B^*|^{1/2})$.
\item [2. ] $R(A^*)\subseteq R(|A+B|^{1/2})$  and $ R(B^*)\subseteq R(|A+B|^{1/2})$.
\een
In this case,  the \textbf{parallel sum} of $A$ and $B$, denoted by 
$A\sump B\in\opunodos$, is : 
\[
\begin{pmatrix}
  A\sump B& 0 \\
  0 & 0
\end{pmatrix}=\left.\begin{pmatrix}
  A & A \\
  A & A+B
\end{pmatrix}\right/ \barr {r}  \\ {(\hiluno\oplus \{0\},\hildos\oplus \{0\})} \earr .
\]
\end{fed}
\begin{rem}\label{comentario1}
Note that the pair $(A,B)$ is weakly summable if and only if the operator matrix $\begin{pmatrix}
  A & A \\
  A & A+B
\end{pmatrix}$ is $(\hiluno\oplus \{0\},\hildos\oplus \{0\})$-weakly complementable. Hence, the parallel sum is well defined.\EOE
\end{rem}

\begin{pro}\label{filmore williams}
Let $A,B\in\opunodos$ be weakly parallel summable operators and let
$E_A$, $E_B$, $F_A$ and $F_B$ be, respectively, the reduced solutions of the equations
\begin{eqnarray}
A\ =&|A^*+B^*|^{1/2}UX \ , \quad &B\ =|A^*+B^*|^{1/2}UX\label{dos},
\\ A^*=&|A+B|^{1/2}X\ , \quad  &B^*=|A+B|^{1/2}X, \label{cuatro}
\end{eqnarray}
where $U$ is the partial isometry of the polar
decomposition of $A+B$. Then:
\begin{equation}\label{filwill}
  A\sump B=F_A^*E_B=F_B^*E_A,
\end{equation}
\end{pro}
\bdem
Note that $|A^*+B^*|^{1/2}U = U|A+B|^{1/2}$. Then, adding 
in (\ref{dos}) and in (\ref{cuatro}), we get
$$
|A+B|^{1/2} =E_A+E_B,    \peso{and} |A^*+B^*|^{1/2}U = F_A^*+F_B^* \ ,
$$
by the uniqueness of the reduced solution. 
By its definition, $A\sump B= A-F_A^*E_A$. Then
\begin{align*}
A\sump B&=A-F^*_AE_A=F^*_A(|A+B|^{1/2}-E_A)=F_A^*E_B \ .
\end{align*}
The other equality follows in a similar way.\edem

\medskip
\begin{cor}
Let $A,B\in\opunodos$ be weakly parallel summable. Then 
$A\sump B=B\sump A$. 
\end{cor}

\medskip

\begin{cor}
Let $A,B\in\opunodos$ be weakly parallel summable and suppose that the
operator $A+B$ has closed range. Then
$
A\sump B= A - A(A+B)^{\dagger}A=A(A+B)^{\dagger}B .
$
\end{cor}

\noi Using Proposition \ref{teorema uno} we
obtain the following analogous result with respect to parallel
sum.

\begin{pro}\label{sucesion para suma paralela.}
Let $A,B\in\opunodos$ be weakly parallel summable and $x\in\hiluno$. Then
there exists a sequence $\{y_n\}$ and  $M >0 $
such that
$$
A(x+y_n)  \xrightarrow[n\rightarrow\infty]{} A\sump B(x) \ , \quad
B(y_n) \xrightarrow[n\rightarrow\infty]{} -A\sump
B(x) \ , 
$$
and  $\pint{|A+B|y_n,y_n} \leq M \ .$
Conversely, if there exist $d \in \hildos\,$, a sequence $\{y_n\}$ in
$\hiluno$  and a real
number $M$ such that
$$
A(x+y_n)  \xrightarrow[n\rightarrow\infty]{} d \ , \quad
B(y_n) \xrightarrow[n\rightarrow\infty]{} -d \ ,
\peso{and}  \pint{|A+B|y_n,y_n} \leq M \ ,
$$
then  $A\sump B(x) = d $. 
\end{pro}

\begin{cor}
Let $A,B\in\opunodos$ be weakly parallel summable. Then
\[
R(A)\cap R(B)\subseteq R(A\sump B)\subseteq \overline{R(A)}\cap
\overline{R(B)}
\]
\end{cor}
\bdem
Given $x\in R(A)\cap R(B)$, let $y,z\in\hiluno$ such that
$Ay=Bz=x$. Then 
$A((y+z)-z)=x  = B(-z)$. 
In consequence, taking $w=y+z$ and $y_n=-z$ for every
$n\in\mathbb{N}$, by Proposition \ref{sucesion para suma
paralela.} we have that $A\sump B(w)=x$, which prove the first
inclusion. The second inclusion follows immediately from Proposition
\ref{sucesion para suma paralela.}. \edem

\subsubsection*{Parallel summable operators}
Let $A,B\in\opunodos$. As we have already pointed out in Remark \ref{comentario1}, 
the  operator pair $(A,B)$ is weakly summable if and only if the block matrix
\[
M = \begin{pmatrix}
  A & A \\
  A & A+B
\end{pmatrix},
\]
is $(\hiluno\oplus\{0\},\hildos\oplus\{0\})$-weakly complementable.
From this point of view, it is natural consider pairs of operators $(A,B)$ such that $M$ 
is $(\hiluno\oplus\{0\},\hildos\oplus\{0\})$-complementable. In this section we shall study such pairs of operators. 

\begin{fed}\rm \label{sumable}
Let $A,B\in\opunodos$. We say that $A$ and $B$ are
\textbf{parallel summable} if
$$
R(A)\subseteq R(A+B)  \peso{and} R(A^*)\subseteq R(A^*+B^{*}) \ .
$$
Note that these conditions imply that
$R(B)\subseteq R(A+B)$ and $R(B^*)\subseteq R(A^*+B^{*})$.
\end{fed}

\begin{rem}
This notion is indeed stronger than weakly summability. For example, take $A,D\in\posop$ 
such that $A\leq D$ but $R(A)\nsubseteq R(D)$. Denote  $B=D-A\in \posop$. 
By Douglas theorem, $R(A)\inc R(A\rai )\subseteq R(D\rai ) = R((A+B)^{1/2})$. 
Similarly, since $B \le D$,  then	also  $R(B)\subseteq R(D\rai ) =R((A+B)^{1/2})$.
However, by hypothesis,  the pair $(A,B)$ can not be parallel summable, 
because $R(A)\nsubseteq R(A+B)=R(D)$. 

\noi  Both notion coincides, for instance, if $R(A+B)$ is closed. In fact, in this case 
$R(A+B)=R(|(A+B)^*|)=R(|(A+B)^*|^{1/2})$ and $R((A+B)^*)=R(|A+B|)=R(|A+B|^{1/2})$.
\EOE
\end{rem}

\noi Clearly, for parallel summable operators, some of the already
proved properties can be improved. Let us mention, for instance, the
following ones.

\begin{pro}\label{sucesion para la strong suma paralela.}
Let $A,B\in\opunodos$ be parallel summable and
$x\in\hiluno$. Then, there exists $y\in\hiluno$ such that
$
A(x+y) = A\sump B(x)$   and  $  By =-A\sump B(x)  . $ 
Moreover, there are projections $P\in L(\hiluno\oplus\hiluno)$, $Q\in L(\hildos\oplus\hildos)$ such that
$R(P^*)=\hiluno\oplus\{0\}$, $R(Q)=\hildos\oplus\{0\}$ and
\[
Q\begin{pmatrix}
  A & A \\
  A & A+B
\end{pmatrix}
=\begin{pmatrix}
  A & A \\
  A & A+B
\end{pmatrix}
P=\begin{pmatrix}
  A\sump B & 0 \\
     0     & 0
\end{pmatrix}.
\]
\end{pro}
\bdem It follows immediately from Proposition \ref{teorema uno para
compatibles}.\edem

\begin{cor}\label{rango de la strong suma paralela}
If $A,B\in\opunodos$ are parallel summable, then $R(A\sump B
)=R(A)\cap R(B)$. 
\end{cor}

\subsubsection*{Parallel substraction}
Given two operators $A,C\in\opunodos$, it seems natural to study the existence of a solution of the equation
$ 
  A\sump X=C,
$ 
that is, if there exists an operator $B\in\opunodos$ parallel summable 
with $A$ such that $A\sump B=C$. For positive
operators this question has been studied, for example, in
\cite{andersonduftrapp}, \cite{PS}, \cite{andersonmortrapp}  and \cite{PEKA}. Clearly, equation $A\sump X=C$ 
may have no
solutions for some pair of operators $(A,C)$. Indeed, 
\coro{rango de la strong suma paralela} implies that, if
equation $A\sump X=C$ 
has a solution, then
$R(C)\subseteq R(A)$ and $R(C^*)\subseteq R(A^*)$, or,
equivalently,  $R(C-A)\subseteq R(A)$ and $R((C-A)^*)\subseteq
R(A^*)$. 

In this section, we shall prove that, if $R(C-A)=R(A)$
and $R((C-A)^*)= R(A^*)$, then there exists a solution of
equation $A\sump X=C$. 
Moreover, we shall
find a distinguished solution, the \textit{parallel substraction of
the operators} $C$ and $A$. 
Given $A\in \opunodos$, let $\cD_A$ be the set of operators
defined by
\[
\cD_A:=\{C\in\opunodos:\;R(C-A)=R(A)\;\;\mbox{and}\;\; R((C-A)^*)=
R(A^*)\}.
\]

\begin{pro}\label{ida y vuelta}
Let $A\in\opunodos$. Then the map $C\mapsto C\sump (-A)$ is a
bijection between the sets $\cD_A$ and $\cD_{-A}$ with
inverse $D\mapsto D\sump A$.
\end{pro}
\bdem
By the definition of summability, it is clear that $-A$ and $C$ are
summable, for every $C \in \cD_A$. 
Let $E$ be the reduced solution of $C-A=AX$ and let $Q$ be a
projection onto $\hildos\oplus\{0\}$ such that
$
Q\begin{pmatrix}
  -A & -A \\
  -A & C-A
\end{pmatrix}= \begin{pmatrix}
  C\sump (-A) & 0 \\
  0 & 0
\end{pmatrix}.
$
Since
$
 \begin{pmatrix}
  C\sump (-A)+A & 0 \\
  0 & 0
\end{pmatrix} = Q\begin{pmatrix}
  0 & -A \\
  -A & C-A
\end{pmatrix}
\ \ \mbox{and}\ \  \begin{pmatrix}
  0 & -A \\
  -A & C-A
\end{pmatrix}\begin{pmatrix}
  -E & 0 \\
  -I & 0
\end{pmatrix}=\begin{pmatrix}
  A & 0 \\
  0 & 0
\end{pmatrix},$ 
we get that 
$$
\begin{pmatrix}
  C\sump (-A)+A & 0 \\
  0 & 0
\end{pmatrix}\begin{pmatrix}
  -E & 0 \\
  -I & 0
\end{pmatrix}= Q  \begin{pmatrix}
  A & 0 \\
  0 & 0
\end{pmatrix} =\begin{pmatrix}
  A & 0 \\
  0 & 0
\end{pmatrix}.
$$ 

\medskip
\noi
This implies that  $R(A)\subseteq R(C\sump (-A)+A)$. As the other
inclusion always holds, we get $R(A)=R(C\sump (-A)+A)$. In a
similar way we can prove that $R(A^*)=R((C\sump (-A)+A)^*)$. Thus,
the mapping $\Phi:\cD_A\to \cD_{-A}$ given by $\Phi(C)=C\sump
(-A)$ is well defined.
To prove that $\Phi^{-1}(D)=D\sump A$, take $C\in\cD_A$
and $x\in\hiluno$. Then there exists $y, z \in\hiluno$ such that
\begin{align*}
C\sump (-A)(x+y)=\big(C\sump (-A)\big)\sump A(x)\ ,\ \  \ \  Ay=-(C\sump (-A))\sump
A(x)\ , 
\end{align*}
\begin{align*}
C(x+y+z)= C\sump (-A)(x+y) \ , \ \ \mbox{and} \ \  Az=-(C\sump
(-A))(x+y).
\end{align*}
So, $A(y+z)=0$ which implies that $C(y+z)=0$. Hence
\[
Cx=C(x+y+z)=C\sump (-A)(x+y)=(C\sump (-A))\sump A(x)\ ,
\]
and the proof is complete.
\edem

\begin{cor}\label{unica solucion}
Let $A\in\opunodos$. For every $C\in\cD_A$, the equation
\[
A\sump X=C
\]
has a solution. Moreover, $C\sump (-A)$ is the unique
solution $X$ which also satisfies 
$$
R(A+X)=R(A) \peso { and } R((A+X)^*)=R(A^*) \ .
$$
\QED
\end{cor}

\begin{fed}\label{definicion de diferencia paralela}\rm
Given $A\in \opunodos$, and $C\in\cD_A$, the \textbf{parallel substraction}
between  the operators $A$ and $C$, denoted by $C\difp A$, is
defined as the unique solution of equation $A\sump X=C$ guaranteed
by Proposition \ref{unica solucion}
\end{fed}

\begin{rem}
Note that, according to our definition, it holds that $C\difp A=C\sump (-A)$; in
particular, several properties of parallel sum are inherited by
parallel substraction.\EOE
\end{rem}

\section{Shorted formulas using parallel sum}
\label{seccion cinco}

\noi In this section we shall prove some formulas for shorted operator using parallel sums and substractions. Throughout this section $\ese$ and $\ete$ will be two fixed closed subspaces of $\hiluno$ and $\hildos$, respectively. The following lemma was proved in \cite{PS} for pairs of positive operators.

\begin{lem}\label{suma paralela con shorted}
Let $A\in\opunodos$ be \stcom. Let $B\in\opunodos$ be such that $(A,B)$  and $(\short{A}{\ese}{\ete},B)$ are parallel sumable and $A\sump B$ is \stcom. Then
\[
\short{A}{\ese}{\ete}\sump
B=\short{\big(A\sump B\big)}{\ese}{\ete} \ .
\]
\end{lem}
\bdem Let $x\in \ese$. By Proposition \ref{sucesion para la strong
suma paralela.}, there exists $y\in\hiluno$ such that
$$
\short{A}{\ese}{\ete}(x+y)= \short{A}{\ese}{\ete}\sump B(x)\
,\peso{and} By=-\short{A}{\ese}{\ete}\sump B(x).
$$
Let $z\in\ese^\bot$ such that
$A(x+y+z)=\short{A}{\ese}{\ete}(x+y)$. Then,
$
A(x+y+z)= \short{A}{\ese}{\ete}\sump B(x)$ and $By=-\short{A}{\ese}{\ete}\sump B(x).
$
So, $A\sump B\begin{pmatrix}   x \\   z \end{pmatrix}=\short{A}{\ese}{\ete}
\sump B(x)$, which implies,  by Proposition \ref{teorema uno},   that 
$\short{(A\sump B)}{\ese}{\ete}x=\short{A}{\ese}{\ete}\sump B(x).$ \edem

\noi The next technical result will be useful throughout this section.

\begin{pro}\label{lema3 para Demetrios}
Let $A\in \opunodos$ be \stcome. If $B\in\opunodos$ satisfies $R(B)=\ete$ 
and $R(B^*)=\ese$, then there exists
$n_0 \in \N$ such that, for every $n\geq n_0$, $A$ and $nB$ are
parallel summable.
\end{pro}

\noi  We need the following lemma.

\begin{lem}\label{lema1 para Demetrios}
Let $A,B\in \opunodos$ such that $R(A)\subseteq\ete$,
$R(A^*)\subseteq\ese$, $R(B)=\ete$ and $R(B^*)=\ese$. 
Then there exists $n_0\in\mathbb{N}$ such that, 
for every $n\geq n_0$, $R(A+nB)=\ete$ and
$R((A+nB)^*)=\ese$.
\end{lem}
\bdem
It suffices to prove that $\ete\subseteq R(A+nB)$ and
$\ese\subseteq R((A+nB)^*)$, because the reverse inclusions hold by
hypothesis. Since $R(B^*) = \ese$, Douglas theorem assures that there exists 
$\alpha>0$ such that $B^*B\geq \alpha P_\ese$.
Then
\begin{align*}
|A+nB|^2&=A^*A+n^2\,B^*B+n\,(A^*B+B^*A) 
\geq 
\Big(\alpha n^2-n\,\|(A^*B+B^*A)\|\Big)P_\ese.
\end{align*}
Take $n_1\in \N$ such that 
$\alpha n^2 > n\,\|(A^*B+B^*A)\|$ 
for every $n\geq n_1\,$.
By Douglas theorem,  $\ese \subseteq
R(|A+nB|)=R((A+nB)^*)$, for $n\geq n_1\,$.
In a similar way, we can prove that there exists $n_2\in\mathbb{N}$
such that for every $n\geq n_2$, $\ete \subseteq R((A+nB))$
holds. Hence, the statement is proved by taking $n_0=\max\{n_1,n_2\}$.\edem

\bdem[Proof of Proposition \ref{lema3 para Demetrios}]
Take, as in Proposition
\ref{teorema uno para compatibles},
 a projection $P \in \opuno $  such that $AP=\short{A}{\ese}{\ete}$ and  $R(P^*)=\ese$. 
 Since $N(B)=\ese^\bot = R(I-P)$, 
it holds that $B(I-P)=0$ and $BP=B$.
By Lemma \ref{lema1 para Demetrios} there exists $n_1\in \N$ such
that  $R((\short{A}{\ese}{\ete}+nB)^*)=\ese$, for every $n\geq n_1$. 
Fix $n \ge n_1$. Given $x\in\hiluno\,$, 
there exists $y\in \ese$ such that
$\short{A}{\ese}{\ete}x= (\short{A}{\ese}{\ete}+nB)y$.
If $z=Py+(I-P)x \in \hiluno\,$, then 
\begin{eqnarray*}
Ax&=&
A\big(Px+(I-P)x\big) = \short{A}{\ese}{\ete}x + A(I-P)x   \\&=&
\big(\short{A}{\ese}{\ete}+nB\big)y+(A+nB)(I-P)x  
\\&=&(A+nB)Py +(A+nB)(I-P)x = (A+nB)z \ .
\end{eqnarray*}
This shows that  $R(A)\subseteq R(A+nB)$. 
Following the same lines, it can be shown that 
there exists $n_2\in \N$ such
that $R(A^*)\subseteq R(A+nB)^*$ for $n \ge n_2\,$. Thus, $A$ and $nB$ are parallel summable for $n\geq \max\{n_1,\ n_2\}$.\edem

\noi Parallel sum may be defined in terms of shorted operators and the next Proposition shows a converse relation.

\begin{pro}\label{Demetrios}
Let $A\in \opunodos$ be \stcom. 
If  $B\in\opunodos$ satisfies 
$R(B)=\ete$ and $R(B^*)=\ese$, then there exists $n_0 \in \N$
such that: 
\ben
\item [\rm 1. ] The pair $(A, \ nB)$ is summable for every $n\geq n_0$, and 
\item [\rm 2. ] $\displaystyle \short{A}{\ese}{\ete}=\lim_{n\to\infty} A\sump (nB) $ (in the norm topology) . 
\een
\end{pro}

\noi Firstly, we shall prove Proposition \ref{Demetrios} in the following particular case.

\begin{lem}\label{lema2 para Demetrios}
Let $A,B\in \opunodos$ be such that $R(A)\subseteq\ete$,
$R(A^*)\subseteq\ese$, $R(B)=\ete$ and $R(B^*)=\ese$. Then,
$A\sump (nB)\convnorm A$ .
\end{lem}
\bdem
Lemma \ref{lema1 para Demetrios} implies that there exists $n_0\geq
1$ such that, for every $n\geq n_0$, $A$ and $nB$ are
parallel summable. Fix $n\geq n_0$. By definition, 
$
A\sump (nB)=A-F_n^*E_n,
$
where $F_n$ and $E_n$ are, respectively, the reduced solution of
$A^*=|A+nB|^{1/2}X$ and $A=|(A+nB)^*|^{1/2}U_nX$, and 
$U_n$ is the partial
isometry of the polar decomposition of $A+nB$.
We shall show that $\|E_n\|\conv 0$ (resp. $\|F_n\|\conv
0$), which clearly implies the desired norm convergence. By Douglas theorem, 
\begin{align}\label{esn}
\|E_n\|&=\inf\big\{\la\in \R \ : \ 
A^*A\leq \la |A+nB|\, \big\} \ \ , \quad n \in \N \ ,
\end{align}
and there exist  $\alpha,\beta>0$ such that $A^*A\leq \beta P_\ese$ and
$B^*B\geq \alpha P_\ese$. Then $(A^*A)^2 \le \beta^2  P_\ese$, and 
\begin{align*}
|A+nB|^2&=A^*A+n^2\,B^*B+n\,(A^*B+B^*A)\\ 
&\geq \big(\alpha n^2-n\,\|(A^*B+B^*A)\|\big)P_\ese 
\geq \frac{\alpha n^2-n\,\|(A^*B+B^*A)\|}{\beta^2} \ (A^*A)^2 \ .
\end{align*}
Recall that L\"owner's theorem states that for every $r\in (0,1]$ $f(x)=x^r$  is operator monotone, i.e. if $0\leq A\leq B$, then $A^r\leq B^r$. Therefore, if $n$ is large enough, 
$$\displaystyle A^*A \le \frac{\beta}{\big(\alpha n^2-n\,\|(A^*B+B^*A)\|\big)\rai} \ \ 
  |A+nB| \ .
$$ 
Hence, \eqref{esn}  implies that
$
\|E_n\| 
\conv \ 0 \ .
$ Analogously, we get that  $\|F_n\|\conv 0\,$. \edem

\bdem[Proof of Proposition \ref{Demetrios}.]
By Proposition \ref{lema3 para Demetrios}, there exists $n_0$ such that for every
$n\geq n_0$ the pairs $(A,nB)$ and $(\short{A}{\ese}{\ete},nB)$
are parallel summable. Since the hypothesis
of Lemma \ref{suma paralela con shorted} are satisfied, for every $n\geq
n_0\,$, it holds that
\[
A\sump nB=\short{(A\sump nB)}{\ese}{\ete}=\short{A}{\ese}{\ete}\sump nB \ .
\]
Then, by Lemma \ref{lema2 para Demetrios} with $\short{A}{\ese}{\ete}$ playing the role of $A$, we get $A\sump
nB\conv\short{A}{\ese}{\ete}$.\edem

\noi Our last result relates parallel sum, parallel substraction and shorted
operators.

\begin{pro}\label{shorted con suma y resta}
Let $A\in\opunodos$ be \stcom, and let $L\in \opunodos$ be such that $R(L)=\ete$ 
and $R(L^*)=\ese$.
Then, there exists $n\in \N$ such that
\begin{enumerate}
  \item [\rm 1.] $A$ and $nL$ are summable.
  \item [\rm 2.] $\short{A}{\ese}{\ete}\in\cD_{-nL}$.
  \item [\rm 3.] $(A\sump nL)\difp nL=\short{A}{\ese}{\ete}$.
\end{enumerate}
\end{pro}
\bdem
The first two assertions follows from Proposition \ref{lema3 para Demetrios} and Lemma \ref {lema1 para Demetrios}, respectively. Since $R\big((A\sump nL)\difp nL\big)\subseteq \ete$ and $R\big(((A\sump nL)\difp nL)^*\big)\subseteq \ese$ then, by Lemma \ref{suma paralela con shorted},
\[
(A\sump nL)\difp nL=\short{\big((A\sump nL)\difp  nL\big)}{\ese}{\ete}=
(\short{A}{\ese}{\ete}\sump nL)\difp nL.
\]
Finally, by Proposition \ref{ida y vuelta},
$
(\short{A}{\ese}{\ete}\sump nL)\difp
nL=(\short{A}{\ese}{\ete}\difp -nL)\sump (-nL)=\short{A}{\ese}{\ete}\ ,
$
and the proof is complete.
\edem

\begin{rem}
Proposition \ref{Demetrios} was proved for positive operators by Anderson and Trapp 
in \cite{andtrapp} and by Pekarev and Smul'jian in \cite{PS}. It was also considered 
by Mitra and Puri who proved  formula 
$\displaystyle \short{A}{\ese}{\ete}=\lim_{n\to\infty} A\sump (nB) $ of \prop {Demetrios} 
for rectangular 
matrices (see \cite{[MP]}). However, their proof can not be extended to infinite 
dimensional Hilbert spaces because it involves generalized inverses which, in our 
setting, only exist for closed range operators. Finally, the reader will find 
a generalization of Proposition \ref{shorted con suma y resta}  for positive 
operators in \cite{PS}.
\end{rem}

\fontsize {9}{8}\selectfont


\begin{thebibliography}{XXXXXX}
\bibitem{anderson} W. N. Anderson Jr., {\em Shorted operators},
 SIAM J. Appl. Math., 20 (1971), 520-525.


\bibitem{andersonduf} W. N. Anderson Jr. and R.J. Duffin, {\em Series
and parallel addition of matrices}, J. Math. Anal. Appl. 11
(1963), 576-594.

\bibitem{andersonduftrapp} W. N. Anderson Jr., R. J. Duffin, and G. E. Trapp, {\em Parallel substraction of matrices}, Proc. Nat. Acad. Sci. 69 (1972), 2530-2531.

\bibitem{andersonmortrapp} W. N. Anderson Jr., D. T. Duffin, and G. E. Trapp, {\em Characterization of parallel substraction}, Proc. Nat. Acad. Sci. 76 (1979), 3599-3601.

\bibitem{andtrapp} W. N. Anderson Jr., G. E. Trapp, {\em Shorted operators
II}, SIAM J. Appl. Math., 28 (1975), 60-71.

\bibitem{Ando} T. Ando, {\em Generalized Schur complements}, Linear
Algebra Appl. 27 (1979), 173-186.

\bibitem{[BG]} A. Ben-Israel and T. N. E. Greville, Generalized inverses.
Theory and applications. Second edition. CMS Books in Mathematics/Ouvrages de
Mathématiques de la SMC, 15. Springer-Verlag, New York, 2003.

\bibitem{[Bo]} R. Bouldin; The product of operators with closed range,
Tohoku Math. J. 25 (1973), 359-363.

\bibitem{butmor} C. A. Butler and T. D. Morley, {\em A note on the shorted operator},
SIAM J. Matrix Anal. Appl., 9 (1988), 147-155.

\bibitem{[CM]} S. L. Campbell and C. D. Meyer Jr., {\em Generalized inverses of linear transformations}, Corrected reprint of the 1979 original. Dover Publications, Inc., New York, 1991. 

\bibitem{[Ca]} D. Carlson, {\em What are Schur complements, anyway?}, Linear
Algebra Appl. 74 (1986), 257-275.

\bibitem{[CaHayn]} D. Carlson and E. V. Haynworth, {\em Complementable and almost definite matrices}, Linear Algebra Appl. 52/53 (1983), 157-176.

\bibitem{CMS2} G. Corach, A. Maestripieri and D. Stojanoff, {\em Schur
complements and oblique projections}, Acta Sci. Math. (Szeged) 67
(2001), 439-459.

\bibitem{CMS1} G. Corach, A. Maestripieri and D. Stojanoff,
{\em Oblique projections  and abstract splines}, J. Approx. Theory 117 (2002),
189--206.

\bibitem{[CMS2]} G. Corach, A. Maestripieri and D. Stojanoff,
{\em Generalized Schur complements and oblique projections}, Linear
Algebra and its  Applications 341 (2002), 259-272.


\bibitem{Co} R. W. Cottle, {\em Manifestations of the Schur complement}, Linear
Algebra Appl. 8 (1974), 189-211.

\bibitem{[De]} F. Deutsch, The angle between subspaces in Hilbert
space, in "Approximation theory, wavelets and applications" (S. P.
Singh, editor), Kluwer, Netherlands, 1995, 107-130.


\bibitem{[Di]} J. Dixmier, Etudes sur les vari\'et\'es et op\'erateurs
de Julia, avec quelques applications, Bull. Soc. Math. France 77
(1949), 11-101.


\bibitem{douglas} R. G. Douglas, {\em On majorization, factorization
 and range inclusion of operators in a Hilbert space}, Proc. Amer. Math. Soc. 17
  (1966), 413-416.

\bibitem{duff} R. J. Duffin, {\em Elementary operations which
generate network matrices}, Proc. Amer. Math. Soc. 11
  (1963), 645-658.

\bibitem{fill} P. A. Fillmore and J. P. Williams, {\em On operator ranges}, Adv. Math.
7 (1971), 254-281.

\bibitem{[Fr]} K. Friedrichs, On certain inequalities and
characteristic value problems for analytic functions and
for functions of two variables. I.
Trans. Amer. Math. Soc. 41 (1937), 321-364 .

\bibitem{Goller} H. Goller, {\em Shorted operators and rank decomposition matrices}, Linear Algebra Appl. 81 (1986), 207-236.

\bibitem {[HJ]}  R. Horn and C. Johnson, {\em Matrix analysis}, Corrected reprint of the 1985 original. Cambridge University Press, Cambridge, 1990. 

\bibitem {[HN]} S. Hassi and K. Nordstr\"om, {\em On projections in a space with an indefinite metric}, Linear Algebra Appl. 208/209 (1994), 401-417.

\bibitem{[Iz]}S. Izumino, The product of operators with closed range
and an extension of the reverse order law, Tohoku Math. J. 34
(1982), 43-52.

\bibitem{[K]} M. G. Krein,  The theory of self-adjoint extensions of
semibounded
Hermitian operators and its applications, Mat. Sb. (N. S.) 20 (62)
(1947), 431-495.



\bibitem{[MS]} P. G. Massey and D. Stojanoff, {\em Generalized Schur complements and 
P-comple\-mentable operators}, Linear Alg. Appl. 393
(2004), 299-318.

\bibitem{[Mitra-Mathew]} T. Mathew and S. K. Mitra, 
{\em Shorted operators and the identification problem -the real case}, IEEE Trans. Circuits and Systems 31 (1984), 299-300.

\bibitem{[Mitra]} S. K. Mitra, 
{\em Shorted operators and the identification problem}, IEEE Trans. Circuits and Systems 29 (1982), 581-583.

\bibitem{[Mitraminus]} S. K. Mitra, 
{\em The minus partial order and the shorted matrix }, Linear Algebra Appl. 83 (1986), 1-27.


\bibitem{[Mitra1]} S. K. Mitra and K. M. Prasad, {\em The nonunique parallel sum}, Linear Algebra Appl. 259 (1997), 77-99.

\bibitem{[Mitra2]} S. K. Mitra and M. L. Puri,{\em On parallel sum and difference of matrices}, J. Math. Anal. Appl. 44 (1973), 92-97.

\bibitem{[MP]} S. K. Mitra and M. L. Puri, {\em Shorted matrices - An
extended concept and some applications}, Linear Alg. Appl. 42
(1982), 57-79.

\bibitem{[MPut]} S. K. Mitra and S. Puntanen, {\em The shorted operator statistically interpreted}, Calcutta Statist. Assoc. Bull. 40 (1990/91), 97-102. 

\bibitem{[MR]} C. R. Rao and S. K. Mitra, {\em Generalized inverses of matrices and its applications}, Wiley, New York, 1971. 

\bibitem{[omelet]} D. V. Ouellette, {\em Schur complements and statistics}, Linear Algebra Appl. 36 (1981), 187-295. 

\bibitem{PEKA} E. L. Pekarev, {\em Shorts of operators and some
extremal problems}, Acta Sci. Math. (Szeged) 56 (1992), 147-163.


\bibitem{PS} E. L. Pekarev and J. L. Smul'jan, {\em Parallel
addition and parallel substraction of operators}, Math. USSR
Izvestija 10 (1976), 351-370.


\bibitem{Smul} J. L. Smul'jan,{\em An operator Hellinger
integral}, Amer. Math. Soc. Translations 22 (1962), 289-337.

\end{thebibliography}
\end{document}